\documentclass[12pt, oneside, a4paper]{article}

\usepackage{amsmath}
\usepackage{amsfonts}
\usepackage{amssymb}
\usepackage{amsthm,mathrsfs}
\usepackage{enumerate}
\usepackage{graphicx}
\newtheorem{theorem}{Theorem}[section]

\theoremstyle{definition}
\newtheorem{remark}[theorem]{Remark}

\title{\textbf{Ramanujan graphs with diameter at most three}}
\author{Mahdi Ebrahimi\footnote{ m.ebrahimi.math@ipm.ir, m.ebrahimi.math@gmail.com}
 \\
 {\small\em  School of Mathematics, Institute for Research in Fundamental Sciences (IPM)},\\{\small\em P.O. Box: 19395--5746, Tehran, Iran},\\{\small\em ORCID ID: 0000-0001-9789-7376}\\
\\
}
\date{}

\begin{document}

\maketitle


\begin{abstract}
For a simple graph $G$, the complement and the line graph of $G$ are denoted by $G^c$ and $L(G)$, respectively.
In this paper, we show that for every simple connected regular graph $G$ with at least $5$ vertices, the graph $\mathcal{R}(G):=L(L(G)^c)^c$ is a Ramanujan graph with diameter at most three.

   \end{abstract}
\noindent {\bf{Keywords:}}  Ramanujan graph, line graph, eigenvalue, diameter. \\
\noindent {\bf AMS Subject Classification Number:}  05C50, 05C75, 05C76.

\section{Introduction}
$\noindent$In this paper, all graphs are assumed to be finite and simple.
 Let $G$ be a graph with vertex set $V(G)=\{\nu_1,\nu_2,\dots, \nu_n\}$ and edge set $E(G)$.
 The complement of $G$ is denoted by $G^c$.
 The \textit{adjacency matrix} of $G$, denoted by $A(G)$, is the $n\times n$ matrix such that the $(i,j)$-entry is $1$ if $\nu_i$ and $\nu_j$ are adjacent, and is $0$ otherwise.
  The \textit{eigenvalues} of $G$ are the eigenvalues of its adjacency matrix $A(G)$.
  The \textit{line graph} $L(G)$ of $G$ is the graph with vertex set $V(L(G)):=E(G)$ in which two vertices are adjacent if and only if their corresponding edges are adjacent, i.e., they are incident to a common vertex in $G$.

 A fundamental subject in graph theory is the study of the spectral gap of a regular graph $G$, that is, the difference between the two largest eigenvalues of $G$ (see \cite{D,A2,B,A1}). Ramanujan graphs are graphs with an optimal spectral gap \cite{A}.
 A \textit{Ramanujan} graph is a connected $k$-regular graph $G$ satisfying $\lambda^*(G)\leq 2\sqrt{k-1}$, where $\lambda^*(G)$ denotes the largest absolute value of an eigenvalue of $G$ distinct from $k$ or $-k$ \cite{13}.
 To construct a Ramanujan graph, some explicit methods are known, for example, we refer the reader to \cite{13,26,21,27,19,6}. These methods are achieved from concepts in linear algebra, number theory, representation theory and the theory of automorphic forms. In this paper, we wish to present a simple combinatorial method to construct a Ramanujan graph.

  \begin{theorem}\label{main}
  For every connected regular graph $G$ with at least $5$ vertices, the graph $\mathcal{R}(G):=L(L(G)^c)^c$ is a Ramanujan graph with diameter at most three.
  \end{theorem}
\begin{remark}
Let $G$ be a connected regular graph with at least $5$ vertices. Set $\mathcal{R}^{(1)}(G):=\mathcal{R}(G)$, and for every positive integer $n\geq 2$, $\mathcal{R}^{(n)}(G):=L(\mathcal{R}^{(n-1)}(G))^c$. Applying Theorem \ref{main}, we deduce that $(\mathcal{R}^{(n)}(G))_{n\in \mathbb{N}}$ is a sequence of Ramanujan graphs with diameter at most three. Therefore the graph $G$ determines a sequence of Ramanujan graphs.
\end{remark}

\section{Ramanujan graphs with diameter at most three}
For a $k$-regular graph $G$ with $n$ vertices and $m$ edges, it is vell-known \cite[Theorem 1.2.16]{spectb1} that
$$P_{L(G)}(x)=(x+2)^{m-n}P_G(x-k+2),$$
where $P_G(x)$ and $P_{L(G)}(x)$ are the characteristic polynomials of $G$ and $L(G)$, respectively. Also we know \cite[Corollary 1.2.14]{spectb1}
$$P_{G^c}(x)=\frac{(-1)^nP_G(-x-1)(x-n+k+1)}{x+1+k}.$$
Now we are ready to prove our main result.\\
\textit{Proof of Theorem \ref{main}:}
Suppose $G$ is a connected $k$-regular graph with $n\geq 5$ vertices. It is clear that $\mathcal{R}(G)$ is a graph with $m:=\frac{1}{8}nk(nk-4k+2)$ vertices. Let $\lambda_1\geq \lambda_2\geq \dots \geq \lambda_m$ be the eigenvalues of $\mathcal{R}(G)$. Then we can see that
$$d:=\lambda_1=\frac{1}{8}(nk-8)(nk-4k+2)+1$$
and
$$\lambda_2=1.$$
It is clear that $1$ is the second largest eigenvalue of $L(G)^c$. Thus $\frac{1}{2}nk-2k$ is the second largest eigenvalue of $L(L(G)^c)$. Therefore
$$\lambda_m=-\frac{1}{2}nk+2k-1.$$
Since $\mathcal{R}(G)$ is a $d$-regular graph and $d>\lambda_2$, using \cite[Theorem 1.2.1]{spectb1} and \cite[Theorem 2.6.2]{S}, we deduce that the graph $\mathcal{R}(G)$ is a connected graph with diameter at most 3. It is obvious that
$$\lambda^*(G)=\mathrm{max}\{|\lambda_i||\,2\leq i\leq m\}=\frac{1}{2}nk-2k+1.$$
Hence
\begin{align}
2\sqrt{d-1}-\lambda^*(G)&=\sqrt{\frac{1}{2}(nk-8)(nk-4k+2)}-\frac{1}{2}(nk-4k+2)\nonumber\\
&=\frac{\frac{1}{2}(nk-4k+2)(\frac{1}{2}nk+2k-9)}{\sqrt{\frac{1}{2}(nk-8)(nk-4k+2)}+\frac{1}{2}(nk-4k+2)}\nonumber\\
&\geq\ 0. \nonumber
\end{align}
Thus $\mathcal{R}(G)$ is a Ramanujan graph and this completes the proof.\qed


\section*{Acknowledgements}
Funding: This research was supported in part
by a grant  from School of Mathematics, Institute for Research in Fundamental Sciences (IPM).


\end{document}